\pdfoutput=1
\RequirePackage{ifpdf}
\ifpdf 
\documentclass[pdftex]{sigma}
\else
\documentclass{sigma}
\fi

\usepackage{mathrsfs}

\def\QQ{\mathbb{Q}}
\def\RR{\mathbb{R}}
\def\CC{\mathbb{C}}

\def\ZZ{\mathbb{Z}}

\def\V{\mathcal{V}}

\def\VV{\mathbb{V}}
\def\nab{\bar{\nabla}}
\def\F{\mathcal{F}}
\def\b{\bullet}
\def\J{\mathcal{J}}
\def\CS{\mathscr{S}}
\def\CX{\mathscr{X}}
\def\CT{\mathscr{T}}
\def\CA{\mathscr{A}}
\def\A{\mathcal{A}}
\def\M{\mathcal{M}}
\def\HH{\mathbb{H}}
\def\H{\mathcal{H}}

\def\w{\omega}
\def\ft{\mathfrak{t}}
\def\fts{\ft^*}
\def\fb{\mathfrak{b}}
\def\fn{\mathfrak{n}}
\def\fp{\mathfrak{p}}
\def\fg{\mathfrak{g}}

\def\l{\lambda}

\def\k12{\mathcal{K}_{\lambda_1,\lambda_2}}
\def\tk12{\tilde{\mathcal{K}}_{\lambda_1,\lambda_2}}
\def\ck12{\check{\mathcal{K}}_{\lambda_1,\lambda_2}}

\def\br{\mathrm}

\numberwithin{equation}{section}

\newtheorem{Theorem}{Theorem}[section]
\newtheorem*{Theorem*}{Theorem}
\newtheorem{Corollary}[Theorem]{Corollary}
\newtheorem{Lemma}[Theorem]{Lemma}
\newtheorem{Proposition}[Theorem]{Proposition}
 { \theoremstyle{definition}
\newtheorem{Definition}[Theorem]{Definition}

\newtheorem{Example}[Theorem]{Example}
\newtheorem{Remark}[Theorem]{Remark} }

\begin{document}
\allowdisplaybreaks

\newcommand{\arXivNumber}{1503.08355}

\renewcommand{\thefootnote}{}

\renewcommand{\PaperNumber}{116}

\FirstPageHeading

\ShortArticleName{Normal Functions over Locally Symmetric Varieties}

\ArticleName{Normal Functions over Locally Symmetric Varieties\footnote{This paper is a~contribution to the Special Issue on Modular Forms and String Theory in honor of Noriko Yui. The full collection is available at \href{http://www.emis.de/journals/SIGMA/modular-forms.html}{http://www.emis.de/journals/SIGMA/modular-forms.html}}}

\Author{Ryan KEAST~$^\dag$ and Matt KERR~$^\ddag$}

\AuthorNameForHeading{R.~Keast and M.~Kerr}

\Address{$^\dag$~Department of Mathematics, University of Toronto, Toronto, Ontario M5S 2E4, Canada}
\EmailD{\href{mailto:ryan.keast@gmail.com}{ryan.keast@gmail.com}}

\Address{$^\ddag$~Department of Mathematics and Statistics, Washington University in St.~Louis,\\
 \hphantom{$^\ddag$}~St.~Louis, MO 63130, USA}
\EmailD{\href{mailto:matkerr@wustl.edu}{matkerr@wustl.edu}}
\URLaddressD{\url{http://www.math.wustl.edu/~matkerr/}}

\ArticleDates{Received May 09, 2018, in final form October 22, 2018; Published online October 26, 2018}

\Abstract{We classify the irreducible Hermitian real variations of Hodge structure admitting an infinitesimal normal function, and draw conclusions for cycle-class maps on families of abelian varieties with a given Mumford--Tate group.}


\Classification{14D07; 14C25; 14M17; 17B45; 32M15; 32G20}

\renewcommand{\thefootnote}{\arabic{footnote}}
\setcounter{footnote}{0}

\vspace{-2mm}

\section{Introduction}\label{S1}

Normal functions are holomorphic horizontal sections of the intermediate Jacobian bundle $J(\mathcal{V})$ associated to a variation of Hodge structure\footnote{Henceforth abbreviated to `VHS' or `variation', with `PVHS' resp.~`PHS' for polarized VHS resp.~HS.} $\mathcal{V}$ of odd weight over a complex analytic manifold. Given a family of smooth projective varieties $\mathscr{X}\overset{\pi}{\to}\mathscr{S}$ and a cycle $\mathfrak{Z}\in {\rm CH}^{p}(\mathscr{X})$ whose restriction $Z_{s}:=\mathfrak{Z}\cdot X_{s}$ to each fiber is homologous to zero, the Abel--Jacobi images $\operatorname{AJ}_{X_{s}}^{p}(Z_{s})$ yield such a section $\nu_{\mathfrak{Z}}$ of $J\big(\mathcal{V}_{\pi}^{2p-1}\big)$.\footnote{Here $\mathcal{V}_{\pi}^{2p-1}$ denotes the VHS with underlying local system $R^{2p-1}\pi_* \QQ$.} These \emph{geometric} normal functions and their singularities are at the heart of the reformulation of the Hodge conjecture by Green and~Griffiths~\cite{GG}, which has brought about a renewal of interest in admissible normal functions and their invariants (cf.~\cite{KP} and references therein).

The most famous example arises from the image $C^{+}$ of an algebraic curve $C$ in its Jacobian. Writing $\imath\colon J(C)\to J(C)$ for the map sending $u\mapsto-u$, and $C^{-}$ for $\imath(C^{+})$, the \emph{Ceresa cycle} $Z_{C}:=C^{+}-C^{-}\in {\rm CH}_{1}(J(C))$ is homologous to zero, with $\operatorname{AJ}_{C}^{g-1}(Z_{C})\in J\big(H_{\rm prim}^{2g-3}(J(C))\big)$ of infinite order for very general~$C$ of genus at least~3~\cite{Ce}. By studying cohomology of the mapping class group, Hain has proved (up to a factor of 2) that the resulting section essentially generates all normal functions over the moduli space $\mathcal{M}_{g}(\ell)$ of curves with level $\ell$ structure (for $\ell\geq3$)~\cite{Ha}.\footnote{Here the underlying VHS $\mathcal{V}\to\mathcal{M}_{g}(\ell)$ is assumed to lie in the Tannakian category generated by $H^{1}(C)$.} By ``essentially'' we mean up to sections of Jacobians of level-1 variations, which are families of abelian varieties. In this paper we shall mostly ignore these, which in the geometric case corresponds to considering only the normal functions that factor through the Griffiths group of cycles modulo algebraic equivalence. We shall also work rationally, i.e., consider torsion normal functions to be zero.

A natural question is what happens over $\mathcal{A}_{g}(\ell)$, or more generally over quotients $\mathcal{A}_{g}(\Gamma)$ of the Siegel domain $D=\br{III}_{g}:=\operatorname{Sp}_{2g}(\mathbb{R})/U(g)$ by a torsion-free congruence subgroup $\Gamma\leq \operatorname{Sp}_{2g}(\mathbb{Q})$.\footnote{$\mathcal{A}_{g}(\ell)$ is the case $\Gamma=\ker (\operatorname{Sp}_{2g}(\mathbb{Z})\to \operatorname{Sp}_{2g}(\mathbb{Z}/\ell\mathbb{Z}))$, which is torsion-free for $\ell\geq3$.} However, Ceresa's normal function does not even survive the (generically $2:1$) map from $\mathcal{M}_{g}(\ell)$ to $\mathcal{A}_{g}(\ell)$, let alone extend to the latter. In fact, a result of Raghunathan~\cite{Ra} implies a~much more general vanishing phenomenon: if~$\mathcal{V}$ is a homogeneous VHS (cf.~Section~\ref{S2.3}) over a locally symmetric variety $\mathrm{X}:=\Gamma\backslash D:=\Gamma\backslash G(\mathbb{R})/K$, associated to a nontrivial $\mathbb{Q}$-irrep\footnote{I.e., irreducible representation defined over $\mathbb{Q}$.}~$V$ of~$G$, and~$G$
has $\mathbb{Q}$-rank $>1$, then $J(\mathcal{V})$ admits no admissible normal functions over any Zariski open subset of~$\mathrm{X}$. Moreover, when the $\mathbb{Q}$-rank is one, the boundary components in the smooth toroidal compactification of~$\mathrm{X}$ are smooth, and so any admissible normal function would have no singularities to study.

This situation improves somewhat when one considers \'etale neighborhoods of $\mathrm{X}$ instead of Zariski ones, motivated by the example of $\mathcal{M}_{3}(\ell)$ minus its hyperelliptic locus (over \smash{$\mathrm{X}=\mathcal{A}_{g}(\ell)$}). To simplify matters, assume that the group~$G$ (of any $\mathbb{Q}$-rank), the Hermitian symmetric domain~$D$, and the representation~$V$ are as above, but also that $G_{\mathbb{R}}$ is simple and $V_{\mathbb{R}}$ irreducible, with $V_{\mathbb{R}}$ of highest weight $\lambda$. (We also consider some cases where $G_{\mathbb{R}}=U(1)\cdot\text{simple}.)$ One would like to classify the pairs $(D,\lambda)$ for which $\mathcal{V}$ has odd weight and the resulting $J(\mathcal{V})\to\mathrm{X}$ admits a normal function after base change to some \'etale neighborhood. For this to happen, a certain cohomology sheaf $\oplus_{j\geq0}\mathcal{H}^{1}(j)$ of \emph{infinitesimal normal functions} on $D$ (cf.\ Section~\ref{S2.1} and Definition~\ref{DefINF}) must not vanish, and it is the pairs with this property that we shall classify
in this paper using a result of Kostant~\cite{Ko} (cf.\ Sections~\ref{S2.2} and~\ref{S2.4}). We should note that this approach has already been carried out by Nori~\cite{No} in the Siegel domain case, so the results below involving $D=\br{III}_{g}$ are not new.

The irreducible Hermitian symmetric domains of noncompact type are recalled in Section~\ref{S2.3} (see the table). We obtain a complete classification of pairs $(D,\lambda)$ admitting an infinitesimal normal function in the case where $D$ is a tube domain (cf.~\cite{Gro})~-- these are the types $\br{I}_{p,p}$, $\br{II}_{2m}$, $\br{III}_n$, $\br{IV}_m$, and $\br{EVII}$ in the table. Ignoring the level $1$ variations, since these all admit an infinitesimal normal function (Example~\ref{Ex2.4}), the calculations in Section~\ref{S3} yield:
\begin{Theorem*}Suppose $D$ is a tube domain, and $\tilde{\mathcal{V}}^{\lambda}$ $($of odd weight at least $3)$ admits an infinitesimal normal function. Then $(D,\lambda)$ is either one of the weight-$3$ Calabi--Yau variations
\begin{gather*}
(\br{I}_{3,3},\omega_{3}), \quad (\br{II}_{6},\omega_{6} ),\quad (\br{III}_{3},\omega_{3} ),\quad (\br{EVII},\omega_{7} )
\end{gather*}
described by Gross~{\rm \cite{Gro}} and Friedman--Laza~{\rm \cite{FL}},\footnote{These are the ``obvious'' members of the list, cf.\ Example~\ref{Ex2.5}. Here ``Calabi--Yau'' (henceforth abbreviated `CY') simply means that the first and last nonzero Hodge numbers are~$1$.} or is a member of one of the infinite families
\begin{gather*}
 (\br{I}_{2,2},\omega_{1}+a\omega_{2}\ \text{or} \ \omega_{3}+a\omega_{2} ),\quad (\br{II}_{4},\omega_{1}+a\omega_{3} \ \text{or} \ \omega_{1}+a\omega_{4} ),\\
 (\br{III}_{1},(2a+1)\omega_{1} ),\quad (\br{III}_{2},\omega_{1}+a\omega_{2} ),\quad (\br{IV}_{2n-1},a\omega_{1}+\omega_{n} ), \quad \text{and}\\
(\br{IV}_{2n-2},a\omega_{1}+\omega_{n-1}\ \text{or} \ a\omega_{1}+\omega_{n} )
\end{gather*}
parametrized by $a\in \mathbb{Z}_{\geq 0}$.
\end{Theorem*}

Our results in the non-tube cases are more subtle and partial; while we enlarge $G$ by a 1-torus and thus allow a variant of the half-twist construction~\cite{vG1} (cf.\ Section~\ref{S4}), we mainly restrict to variations of CY type or appearing in the cohomology of abelian varieties of generalized Weil type (cf.\ Definition~\ref{DefGW}).

Returning to $\mathcal{A}_{g}(\ell)$ (or more generally $\mathcal{A}_{g}(\Gamma)$), we can ask once more about the loci ``suppor\-ting'' higher-weight normal functions.\footnote{The question is independent of $\ell$ (or $\Gamma$), but we need $\Gamma$ torsion-free to be able to speak of a homogeneous VHS over \emph{all} of $\mathcal{A}_{g}(\Gamma)$.} By this we mean that there is an irreducible \emph{homogeneous} variation~$\mathcal{V}$ of odd level $>1$ over $\mathcal{A}_{g}(\Gamma)$, a subvariety $\mathscr{S}\overset{\imath}{\hookrightarrow}\mathcal{A}_{g}(\Gamma)$ with \'etale neighborhood $\mathscr{T}\overset{\jmath}{\to}\mathscr{S}$, and a nontorsion horizontal section of $J(\jmath^{*}\imath^{*}\mathcal{V})$ over $\mathscr{T}$. Alternatively, one could consider only those~$\mathcal{V}$ appearing in the cohomology of the universal abelian variety. In either case, the loci $\overline{\imath(\mathscr{S})}$ are proper subvarieties for $g>3$~\cite{No}, and are expected to be more like $\mathcal{M}_{g}(\Gamma)$ than locally symmetric subvarieties. Indeed, $\mathcal{M}_{g}(\Gamma)$ supports in this sense the Ceresa normal function for any $g$; another example is the Fano normal function supported by the locus in~$\mathcal{A}_{5}(\Gamma)$ of intermediate Jacobians of cubic threefolds~\cite{CNP}. The approach in this article is consequently analogous to the search for special subvarieties of the Torelli locus $\overline{\mathcal{M}_{g}}\subset\mathcal{A}_{g}$~\cite{MO}, which are expected not to exist for large enough~$g$. So the following vanishing result, which follows from our classification, should not be surprising:
\begin{Theorem*}[Theorems~\ref{thm1} and~\ref{thm2}] The loci of Weil resp.\ quaternionic abelian varieties\footnote{See Definitions~\ref{DefQAV} and~\ref{DefWAV} for quaternionic resp.\ Weil abelian varieties, and~\eqref{eq 3.1} for the reduced Abel--Jacobi map.} \linebreak in~$\mathcal{A}_{g}(\Gamma)$ do not support higher weight normal functions for $g>6$ resp.~$8$. In particular, the image of the reduced Abel--Jacobi map
\begin{gather*}
\overline{\operatorname{AJ}}^{r}\colon \ \operatorname{Griff}^{r}(A)_{\mathbb{Q}}\to\overline{J}^{r}(A)_{\mathbb{Q}}
\end{gather*}
 is zero for a very general Weil resp.\ quaternionic abelian variety $A$ of dimension at least $8$ resp.~$10$.
\end{Theorem*}
There are related results for cycles in certain codimensions for the generalized Weil abelian varieties over $\br{I}_{p,q}$ ($p+q=n+1$, $p<q$) studied in Section~\ref{S4.1}, which get weaker as $|q-p|$ grows (cf.\ Theorem~\ref{thm3}). However, both here and for the domains $\br{IV}_{m}$ parametrizing ``spin abelian varieties'' (for instance, those arising from the Kuga--Satake construction), one has infinitesimal normal functions
for arbitrarily large~$n$ and~$m$.

For each of the ``infinitesimal normal functions'' described in this paper, the next immediate problem is to determine whether it comes from an actual (admissible) normal function, and if so, to geometrically realize it and determine\footnote{Outside the exceptional cases $(\br{EVII},\omega_7)$ and $\big(\br{EIII},\omega_1 \big\{ \tfrac{1}{6} \big\} \big)$, which do not embed in~$\mathcal{A}_g$.} its locus of support in $\mathcal{A}_{g}(\Gamma)$. As an added incentive, whenever the geometric realization is as a cycle on some finite pullback of the universal abelian family over a locally symmetric subvariety $D(\Gamma')\subseteq\mathcal{A}_{g}(\Gamma)$, one obtains as a corollary the infinite generation of the relevant Griffiths group of a very general member of this family. This was observed by Nori~\cite{No2} for type $\br{III}_{3}$ (i.e., $\mathcal{A}_{3}(\Gamma)$ itself) using the Ceresa cycle; the general argument is practically identical to that in [op.\ cit.], but with~\cite{Ra} replacing the use of the congruence subgroup theorem.

We briefly carry out the geometric realization for $(D,\lambda)= (\br{I}_{3,3},\omega_{3} )$ in Section~\ref{S3.5}, using the fact that for certain families of Weil abelian $6$-folds, the general member is a Prym variety associated to an unramified $3:1$ curve covering with base curve of genus $g=4$~\cite{F,Sc}. See Theorem~\ref{thm Griff Weil} for the application to Griffiths groups. The locus of support of the resulting normal function is (at least) the $9$-dimensional locally symmetric quotient of~$D$ in~$\mathcal{A}_{6}(\Gamma)$. The construction works for $g>4$ as well, but a dimension count shows that the locus of support is a \emph{proper} subvariety of the relevant $\Gamma\backslash\br{I}_{p,p}$.

With an admissible normal function in hand, one can also try to determine its zero locus, and its singularities along intersections of toroidal boundary components. But even in the absence of this, one can probably still use Kostant's result, replacing $\mathfrak{g}$ by subalgebras of the form $\mathfrak{sl}_{2}^{\oplus a}$, to compute spaces of singularity classes for normal functions over homogeneous variations.

An equally intriguing prospect is to try to extend the computation of infinitesimal normal function spaces (or singularity classes) to the nonclassical case, which is more relevant to the Green--Griffiths program~\cite{GG}. Here the homogeneous families of Hodge structures only become variations upon restriction to the image of a period map. However, the computation only requires an infinitesimal variation as input, so it may already be interesting to start with tangent spaces to the Schubert VHS of~\cite{Ro}.

\section{Infinitesimal normal functions and Kostant's theorem} \label{S2}

This section contains the abstract underpinnings of the calculations in Sections~\ref{S3} and~\ref{S4}. Apart from Section~\ref{S2.4}, it is mostly expository and contains the material on infinitesimal invariants, Lie algebra cohomology, and Hermitian VHS that will be used subsequently.

\subsection{Normal functions and infinitesimal invariants} \label{S2.1}

Let $\mathcal{V}\to\mathscr{S}$ be a polarized $\QQ$-VHS over a~complex manifold, pure of weight~$-1$. By abuse of notation, $\mathcal{V}$ will denote the underlying holomorphic vector bundle and its sheaf of sections, with Hodge filtration~$\F^{\b}$; $\VV$ is the underlying $\QQ$-local system, $Q\colon \VV\times\VV\to\QQ_{\CS}$ the polarization, and~$\nabla$ the Gauss--Manin connection.

Define a filtered complex of sheaves\begin{gather*} 
C^{\b}:= ( \V\otimes \Omega^{\b}_{\CS},\nabla ) ,\qquad F^p C^{\b}:= ( \F^{p-\b}\V\otimes \Omega^{\b}_{\CS},\nabla )
\end{gather*}on $\CS$, so that the differential of each
\begin{gather*}
{\rm Gr}_{F}^{p}C^{\b}=\big({\rm Gr}_{\F}^{p-\b}\V\otimes\Omega_{\CS}^{\b},\nab\big)
\end{gather*}
is $\mathcal{O}_{\CS}$-linear. Consider the exact sequence\begin{gather} \label{eq 1.2}
0\to F^{0}C^{\b}\to C^{\b}/\VV\to C^{\b}/\big(F^{0}C^{\b}+\VV\big)\to0 ;
\end{gather}the sheaf of (quasi-)horizontal sections of the Jacobian bundle $J(\V)$ is
\begin{gather*}
\J_{\rm hor}^{\QQ}=\H^{0}_{\nabla}\big(C^{\b}/\big(F^{0}C^{\b}+\VV\big)\big).
\end{gather*}
The space of normal functions associated to $\V$ is
\begin{gather*}
\operatorname{NF}_{\CS}(\V):=\Gamma\big(\CS,\J_{\rm hor}^{\QQ}\big),
\end{gather*}
with infinitesimal invariant
\begin{gather*}
\delta\colon \ \operatorname{NF}_{\CS}(\V)=\HH^{0}\left(\CS,\frac{C^{\b}}{F^{0}C^{\b}+\VV}\right)\to\HH^{1}\big(\CS,F^{0}C^{\b}\big)\to\Gamma\big(\CS,\H_{\nabla}^{1}\big(F^{0}C^{\b}\big)\big)
\end{gather*}
arising from~\eqref{eq 1.2}.
\begin{Lemma} If $\mathcal{H}_{\nabla}^{0}\big(F^{0}C^{\b}\big)=0$, then $\delta$ is injective.\end{Lemma}
\begin{proof}The assumption implies that we only need to check injectivity into $\mathbb{H}^{1}\big(\CS,F^{0}C^{\b}\big)$. By the theorem of the fixed part~\cite{S}, the assumption also implies that $H^{0}(\CS,\mathbb{V})=\{0\}$. Since $\HH^{0}(\CS,C^{\b}/\VV)=H^{0}(\CS,\VV_{\CC}/\VV)=H^{0}(\CS,\VV)\otimes\tfrac{\CC}{\QQ}$,
the result now follows from the long-exact sequence of \eqref{eq 1.2}.\end{proof}

\begin{Lemma}If $\H_{\nabla}^{1}\big(F^{0}C^{\b}\big)$ and $\H_{\nabla}^{0}\big(F^{0}C^{\b}\big)$ vanish, then $\operatorname{NF}_{\CT}(\jmath^{*}\V)=\{0\}$ for any \'etale neighborhood $\jmath\colon \CT\to\CS$.\end{Lemma}
\begin{proof}Vanishing of $\H^{0}\big(F^{0}C^{\b}\big)=\VV_{\CC}\cap\F^{0}$ implies (again by \cite{S}) that $H^{0}(\CT,\jmath^{*}\VV)=\{0\}$; now apply the previous lemma.
\end{proof}

Now writing $\H^{k}(j):=\H_{\nab}^{k}\big({\rm Gr}_{F}^{j}C^{\b}\big)$, consider the spectral sequence\begin{gather} \label{eq 2.2}
\mathcal{E}_{1}^{q,p}:= \begin{cases}\H^{p+q}(p), & p\geq0, \\0, & p<0\end{cases}
\end{gather}converging to $\H_{\nabla}^{*}\big(F^{0}C^{\b}\big)$. We recover at once the result of Green and Voisin \cite{Gre}, in the form stated by Nori~\cite{No}:
\begin{Proposition}\label{Nori}If all the $\H^{0}(j)$ and $\H^{1}(j)$ vanish for $j\geq0$, then $\operatorname{NF}_{\CT}(\jmath^{*}\V)=\{0\}$ for any \'etale neighborhood $\jmath\colon \CT\to\CS$.
\end{Proposition}

Two obvious situations in which the vanishing conditions \emph{fail} are those of VHS of level one, or level three and ``CY type''. Recall that if
\begin{gather*}
\V=\F^{p_{\min}}\V\supsetneq\F^{p_{\min}+1}\V\supsetneq\cdots\supsetneq\F^{p_{\max}}\V\supsetneq\F^{p_{\max}+1}\V=\{0\}
\end{gather*}
then the level $\ell(\V)=p_{\max}-p_{\min}$; for weight $-1$ we have $\ell(\V)=2p_{\max}+1=-2p_{\min}-1$. Write $h^{p}:=\dim({\rm Gr}_{\F}^{p}\V)$ for the Hodge numbers, and $\underline{h}=\big(h^{p_{\min}},\ldots,h^{p_{\max}}\big)$; put $d:=\dim_{\CC}\CS$. By naively computing ranks we have
\begin{Example}[level 1] \label{Ex2.4}If $\underline{h}=(n,n)$ and $d>1$, then $\H^{1}(0)\neq\{0\}$. If $d=1,$ then $\H^{1}(1)\neq\{0\}$.
\end{Example}

\begin{Example}[level 3 CY] \label{Ex2.5} If $\underline{h}=(1,n,n,1)$ and $1<d<2n$, then $\H^{1}(0)\neq\{0\}$. If $d=1$, then $\H^{1}(2)\neq\{0\}$. (All CY VHS we consider have $d\leq n$.)
\end{Example}

This makes the broad vanishing results we obtain in this paper somewhat surprising. It also explains why we have to ignore the level-one cases below.

Finally, given a smooth family of varieties $\pi\colon \CX\to\CS$, let $Z\in Z^{r}(X_{s_{0}})_{\QQ,{\rm hom}}$ be a homologically trivial cycle on a very general fiber. Suppose its class in the Griffiths group
$\operatorname{Griff}^{r}(X_{s_{0}})=\tfrac{Z^{r}(X_{s_{0}})_{\QQ,{\rm hom}}}{Z^{r}(X_{s_{0}})_{\QQ,{\rm alg}}}$ has nonzero image under the \emph{reduced Abel--Jacobi map} \begin{gather} \label{eq 3.1}
\overline{\operatorname{AJ}}^r_{X_{s_0}}\colon \ \operatorname{Griff}^r (X_{s_0} ) \to \overline{J}^r(X_{s_0}):=\operatorname{Ext}^1_{{\rm MHS}}\big( \QQ(0),\overline{H}^{2r-1}(X_{s_0},\QQ(r))\big) ,
\end{gather}where $\overline{H}^{2r-1}(X_{s_{0}})$ is the quotient of $H^{2r-1}(X_{s_{0}})$ by its maximal level-one sub-HS. Spreading out yields a cycle $\mathfrak{Z}\in Z^{r}(\CX_{\CT})_{\QQ}$ on an \'etale pullback $\CX_{\CT}\overset{\pi_{\CT}}{\to}\CT$ meeting fibers properly (with each $\mathfrak{Z}\cdot X_{t}\underset{\rm hom}{\equiv}0$), and a normal function
\begin{gather*}
\nu_{\mathfrak{Z}}\in\operatorname{NF}_{\CT}(\V),\qquad \V=\overline{R}^{2r-1}(\pi_{\CT})_{*}\QQ_{\CX_{\CT}}(r)\otimes\mathcal{O}_{\CT},
\end{gather*}
where $\overline{R}$ denotes the quotient by the maximal level-one sub-VHS. So we have the
\begin{Corollary}If all the $\H^{0}(j)$ and $\H^{1}(j)$ vanish for $j\geq0$, then $\overline{\operatorname{AJ}}_{X_{s_{0}}}^{r}(Z)=0$ in~\eqref{eq 3.1}.
\end{Corollary}
While there is no converse result, nonvanishing of $\H^{1}(0)$ in particular seems to be a good predictor of the existence of interesting cycles. Note that its nonvanishing for level-one VHS~$\V$ has a geometric ``origin''. Namely, these VHS correspond to the $H^{1}$ (or~$H^{2D-1}$) of families of abelian $D$-folds $\CA\overset{\pi}{\to}\CS$. The existence of nontorsion points on the geometric generic fiber over~$\overline{\CC(\CS)}$ (the algebraic closure of the function field) yields nontrivial geometric normal functions in $\operatorname{NF}_{\CT}(\jmath^{*}\V$) over some \'etale neighborhood~$\CT$. This explains why we only consider the Griffiths group.

\subsection{$\mathfrak{n}$-cohomology of finite-dimensional representations} \label{S2.2}

Let $\fg$ be a complex semisimple Lie algebra of rank $n$, $\fb\supset\ft$ Borel and Cartan subalgebras, $\Delta=\Delta(\fg,\ft)\subset\fts$ the corresponding roots. Denote by $\Delta^{+}=\Delta(\fb)$ the positive roots, $\Sigma=\{\sigma_{1},\ldots,\sigma_{n}\}\subset\Delta^{+}$ the simple roots, $\Omega=\{\w_{1},\ldots,\w_{n}\}\subset\fts$ the fundamental weights, and $\Lambda$ the (weight) lattice they generate. The Killing form $B(X,Y)=\operatorname{Tr} (\operatorname{ad}X\circ\operatorname{ad}Y )$ on $\fg$ induces a~symmetric bilinear form $\langle\,,\,\rangle$ on~$\Lambda$, a~particular orthonormal basis of which (as in \cite[Appendix~C]{Kn}) will be denoted by $\{e_{i}\}$. We have in particular $\langle\w_{i},\sigma_{j}\rangle=\frac{1}{2}\langle\sigma_{j},\sigma_{j}\rangle\delta_{ij}$.

By the theorem of the highest weight, the irreducible representations $\big\{V^{\lambda}\big\}$ of $\fg$ (of finite dimension) are parametrized by their highest weight $\lambda$; there is a 1-to-1 correspondence between irreps and weights of the form $\lambda=\sum\limits_{i=1}^{n}m_{i}\w_{i}$ with $m_{i}\geq0$. Fix an element ${\tt E}\in\ft$ with all $\frac{1}{2}{\tt E}(\sigma_{i})$ non-positive and integral, and let $\fg=\oplus_{j\in\ZZ}\fg^{j,-j}$ be the decomposition into $\operatorname{ad}({\tt E})$-eigenspaces with eigenvalue $2j$. For any representation $(V,\rho)$, there is a corresponding decomposition into $\rho({\tt E)}$-eigenspaces, which can have odd or fractional eigenvalues. Write $\fp=\oplus_{j\geq0}\fg^{j,-j}$, $\fn=\oplus_{j<0}\fg^{j,-j}$ (so that $\Delta(\fn)\subset\Delta^{+}$), and $\fg^{0}=\fg^{0,0}$. We also denote $\Delta_{0}=\Delta\big(\fg^{0},\ft\big)$, $\Delta_{0}^{+}=\Delta_{0}\cap\Delta^{+}$, and $\big\{V_{0}^{\xi}\big\}$ for the highest-weight irreps of~$\fg^{0}$.

Our main computational tool will be a result of Kostant (cf.~\cite[Theorem~5.14]{Ko}). It computes the decomposition of the cohomologies $H^k(\mathfrak{n},V^{\lambda})$ of the natural complex
\begin{gather*}
0\to V^{\lambda}\to\fn^{\vee}\otimes V^{\lambda}\to\wedge^{2}\fn^{\vee}\otimes V^{\lambda}\to\cdots
\end{gather*}
under the action of $\fg^{0}$. To state a version of it (Proposition~\ref{Kostant} below), let $W=W(\fg,\ft)$ and $W_{0}=W\big(\fg^{0},\ft\big)$ be the Weyl groups, and consider the set
\begin{gather*}
W^{0}:=\big\{ w\in W\,|\,w(\Delta^{+})\supseteq\Delta_{0}^{+}\big\}
\end{gather*}
of minimal-length representatives of cosets $W_{0}\backslash W$ (with length $|w(\Delta^{+})\cap\Delta^{-}|$), which we partition into\begin{gather*} 
W^0(j) = \big\{ w\in W^0 \,|\, \text{length}(w)=j \big\} .
\end{gather*}Finally, write $\rho:=\tfrac{1}{2}\sum\limits_{\delta\in\Delta^{+}}\delta=\sum\limits_{i=1}^{n}\w_{i}$,
and\begin{gather*} 
w\cdot\lambda:=w(\lambda+\rho)-\rho .
\end{gather*}
\begin{Proposition}\emph{\label{Kostant}} As a $\fg^{0}$-module, $H^{k}\big(\fn,V^{\lambda}\big)\cong\oplus_{w\in W^{0}(k)}V_{0}^{w\cdot\lambda}$.
\end{Proposition}
Now suppose $\fg\supset\ft$ is the complexification of $\fg_{\RR}\supset\ft_{\RR}$, with $\ft_{\RR}$ \emph{compact} (and $\fg_{\RR}$ not). Writing $\Delta=\Delta_{c}\amalg\Delta_{n}$ for the decomposition into compact and noncompact roots, we assume $\tfrac{1}{2}{\tt E}(\Delta_{c})\subset2\ZZ$ and $\tfrac{1}{2}{\tt E}(\Delta_{n})\subset2\ZZ+1$. (This ensures polarizability of the Hodge structures induced on $\tilde{V}$ and~$\fg$, cf.~\cite[pp.~90--93]{GGK}.) Let~$\tilde{V}$ be an irreducible representation (of finite dimension). We then have either $\tilde{V}_{\CC}\cong V^{\lambda}$ (\emph{real} case) or $\tilde{V}_{\CC}\cong V^{\lambda}\oplus\overline{V^{\lambda}}$ as $\fg$-modules. Let $w_{0}\in W$ denote the unique element with $w_{0}(\Delta^{+})=\Delta^{-}$, and write $-\tau$ for the induced involution on $\Lambda$; for $\fg$ simple of Cartan type other than $A_{n}$, $D_{2m+1}$, or $E_{6}$, $\tau={\rm id}_{\Lambda}$. In the non-real case, we have $\overline{V^{\lambda}}\cong V^{\tau(\lambda)}$, and we say $V^{\lambda}$ is \emph{complex} resp.\ \emph{quaternionic}
when $\lambda\neq\tau(\lambda)$ resp.\ $\lambda=\tau(\lambda)$. To distinguish the real and quaternionic cases, we use the following
\begin{Proposition}[\cite{GGK}] For $\lambda=\sum M_{i}\sigma_{i}$ with $\tau(\lambda)=\lambda$, $V^{\lambda}$ is real iff $\sum\limits_{\sigma_{i}\, {\rm compact}}M_{i}\in\ZZ$.
\end{Proposition}
Given $\l=\sum m_{i}\w_{i}$, $m_{i}\geq0$, we shall write\begin{gather*} 
\tilde{V}^{\l} := \begin{cases}V^{\l}, & \text{real case}, \\ V^{\l}\oplus V^{\tau(\l )}, & \text{complex/quaternionic cases.}\end{cases}
\end{gather*}It is to these representations that we shall apply Proposition~\ref{Kostant}. Note that since $\tilde{V}^{\l}$ has an underlying $\RR$-vector space $\tilde{V}_{\RR}^{\l}$ (which is irreducible as a representation
of $\fg_{\RR}$), it can be viewed (up to Tate twist) as an $\RR$-Hodge structure via the action of ${\tt E}$, whose eigenvalues are regarded as ``$p-q$'' (on Hodge type~$(p,q)$).

Since differences of weights of $V^{\l}$ lie in the root lattice (which ${\tt E}$ sends into $2\ZZ$) and $-\tau(\l)=w_{0}(\l)$ is a weight of $V^{\l}$, $\tilde{V}^{\l}$ \emph{can be placed in weight~$-1$ iff ${\tt E}(\l)$ is odd}. In this case\begin{gather*} 
\big( \tilde{V}^{\l}\big) ^{p,-p-1} :=\big\{ (2p+1)\text{-eigenspace of }{\tt{E}}\text{ on }\tilde{V}^{\l}\big\} ,
\end{gather*}and (using $\overline{V^{\l}}\cong V^{\tau(\l)}$) the level is\begin{gather*} 
\ell\big( \tilde{V}^{\l}\big) = \max \{-\tt{E}(\l),-\tt{E}(\tau(\l)) \} .
\end{gather*}Moreover, when ${\tt E}(\l)$ is odd there exists (up to scale) a~unique $\fg$-invariant alternating bilinear form $Q$ on $\tilde{V}_{\RR}^{\l}$, which polarizes this Hodge structure. Note that in the non-real cases, $Q$~pairs~$V^{\l}$ and~$V^{\tau(\l)}$.
\begin{Remark} The level of the complex summands (as complex Hodge structures) is\begin{gather} \label{eq 6.3}
\ell \big(V^{\l}\big) \big( =\ell\big(V^{\tau(\l)}\big) \big) = -\tfrac{1}{2} ( \tt{E}(\l)+\tt{E}(\tau(\l)) ) .
\end{gather}This is the minimal level that can be achieved using half-twists, which are discussed in Section~\ref{S4}.
\end{Remark}

\subsection{Homogeneous VHS over locally symmetric varieties} \label{S2.3}

Let $G$ be a semisimple $\QQ$-algebraic group of Hermitian type, such that $G_{\RR}$ has a compact maximal torus $T_{\RR}$; and write $\fg_{\RR}\supset\ft_{\RR}$ for the Lie algebras (with complexifications
$\fg\supseteq\ft$). Choose a cocharacter $\chi_{0}\colon \mathbb{G}_{m}\to T_{\CC}$ such that ${\tt E}:=\chi_{0}'(1)$ satisfies\footnote{One should also assume ${\tt E}$ doesn't project to zero in any simple
factor of~$\fg$, but in the calculations below $\fg$ will be simple.}\begin{gather} \label{eq 7.1}
{\tt{E}}(\Delta_c) = 0,\qquad {\tt{E}}(\Delta_n)=\{\pm 2\} .
\end{gather}Denoting the composition $\mathbb{G}_{m}\overset{\chi_{0}}{\to}T_{\CC}\hookrightarrow G_{\CC}$ by~$\varphi_{0}$, the corresponding \emph{Hermitian symmetric domain} is the orbit under conjugation by the group of real points~$G(\RR)$:
\begin{gather*}
D:=G(\RR).\varphi_{0}\cong G(\RR)/G^{0}(\RR),
\end{gather*}
with $G^0(\RR)\leq G(\RR)$ a maximal compact subgroup. Taking $\Gamma\leq G(\QQ)$ a torsion-free congruence subgroup, the \emph{locally symmetric variety} $\mathrm{X}:=\Gamma\backslash D$ is in fact a~quasi-projective variety by the Baily--Borel theorem.

Recall that $B$ denotes the Killing form on $\mathfrak{g}$ (Section~\ref{S2.2}). The $G(\RR)$-orbit of $\operatorname{Ad}\circ\varphi_{0}\colon \mathbb{G}_{m}\to\operatorname{Aut}(\fg_{\CC},-B)$ gives a $(-B)$-polarized $\QQ$-VHS of weight zero and level two over $D$, with Hodge decompositions $\fg=\oplus_{j=-1,0,1}\fg_{\operatorname{Ad}(g){\tt E}}^{j,-j}$ at $g\varphi_{0}g^{-1}\in D$; this descends to a VHS over~$\mathrm{X}$. More
generally (with $\mathbb{A}=\QQ$ or $\RR$), given an $\mathbb{A}$-representation $\rho\colon G_{\mathbb{A}}\to \operatorname{Aut}(V_{\mathbb{A}},Q)$ ($Q$ a $(-1)^{k}$-symmetric bilinear form) such that $\rho\circ\varphi$ is an $\mathbb{A}$-PHS, we get a homogeneous $\mathbb{A}$-PVHS $\V_{\mathbb{A}}$ (with weight of parity~$k$) over~$\mathrm{X}$, called a \emph{Hermitian VHS}. In every case the Hodge decomposition is induced by $d\rho(\operatorname{Ad}(g){\tt E})$; that is, a weight subspace $V_{\xi}\subset V^{p,q}$ $\iff$ ${\tt E}(\xi)=p-q$.

For groups other than $E_{6}$ and $E_{7}$, most\footnote{See Section~\ref{S3.4} for discussion of one possible exception.} Hermitian $\QQ$-PVHS arise from the relative cohomology of various canonical families of abelian varieties; in some cases (cf.~\cite{FL,R}) one has also families of CY varieties. While these structures may matter for questions of geometric origin (i.e., algebraic cycles), the behavior of the $\big\{\H^{k}(j)\big\}$ depends only on~$\V_{\RR}$. Moreover, any Hermitian PVHS of weight~$-1$ on $\mathrm{X}$ decomposes \emph{over $\RR$} into a direct sum of the irreducible homogeneous real variations $\tilde{\V}_{\RR}^{\l}$ arising (as above) from the
PHS on $\big(\tilde{V}^{\l},\tilde{\rho}^{\l}\big)$ described at the end of Section~\ref{S2.2}. Indeed, since $\Gamma\leq G(\RR)$ is Zariski-dense, we have an equivalence between real Hermitian VHS $\V_{\RR}$ over~$\mathrm{X}$ and representations of~$G_{\RR}$.

Before turning to our main calculation, we remind the reader what forms $D$ and ${\tt E}$ can take when $G_{\CC}$ is simple. First, \eqref{eq 7.1} forces $\Delta_{n}\cap\Sigma$ to be a singleton $\{\sigma_{{\tt I}}\}$, which must additionally be a \emph{special} simple root. That is, if we write $\fg=V^{\l_{\operatorname{ad}}}$, $\l_{\operatorname{ad}}=\sum\limits_{i=1}^{n}M_{i}^{\operatorname{ad}}\sigma_{i}$, then $\sigma_{{\tt I}}$ must be one of the simple roots $\sigma_{i}$ for which $M_{i}^{\operatorname{ad}}=1$. Further, the choice of $\sigma_{{\tt I}}$ determines: the decomposition $\Delta=\Delta_{c}\amalg\Delta_{n}$, and thus the real form $\fg_{\RR}$; and the Hodge structure at the base point $\varphi_{0}$, by\begin{gather*} 
{\tt{E}}(\sigma_{\tt{I}})=-2 ,\qquad {\tt{E}}(\sigma_j)=0\qquad \forall\, j\neq\tt{I} .
\end{gather*}In this way, the isomorphism classes of the irreducible Hermitian symmetric domains of noncompact type are parametrized by the choice of Cartan type and special simple root:\begin{gather*}\label{table}
\begin{tabular}{|c|c|c|c|c|}
\hline
$D$ & $(R,\sigma_{{\tt I}})$ & $G_{\RR}$ & $d$ & range\tabularnewline
\hline
\hline
$\br{I}_{p,n-p+1}$ & $(A_{n},\sigma_{p})$ & ${\rm SU}(p,n-p+1)$ & $p(n-p+1)$ & $\begin{matrix}
n\geq2\\
1\leq p\leq\big\lfloor \tfrac{n+1}{2}\big\rfloor
\end{matrix}$\tsep{3pt}\bsep{7pt}\tabularnewline
\hline
$\br{II}_{n}$ & $(D_{n},\sigma_{n})$ & $\operatorname{Spin}^{*}(2n)$ & $\tfrac{1}{2}n(n-1)$ & $n\geq4$\tsep{3pt}\bsep{3pt} \tabularnewline
\hline
$\br{III}_{n}$ & $(C_{n},\sigma_{n})$ & $\operatorname{Sp}(2n,\RR)$ & $\tfrac{1}{2}n(n+1)$ & $n\geq1$\tsep{3pt}\bsep{3pt}\tabularnewline
\hline
$\br{IV}_{2n-1}$ & $(B_{n},\sigma_{1})$ & $\operatorname{Spin}(2,2n-1)$ & $2n-1$ & $n\geq3$\tsep{1pt}\bsep{1pt}\tabularnewline
\hline
$\br{IV}_{2n-2}$ & $(D_{n},\sigma_{1})$ & $\operatorname{Spin}(2,2n-2)$ & $2n-2$ & $n\geq4$\tsep{1pt}\bsep{1pt}\tabularnewline
\hline
$\br{EIII}$ & $(E_{6},\sigma_{1})$ & $E_{6}(-14)$ & $16$ & -----\tsep{1pt}\bsep{1pt}\tabularnewline
\hline
$\br{EVII}$ & $(E_{7},\sigma_{7})$ & $E_{7}(-25)$ & $27$ & -----\tsep{1pt}\bsep{1pt} \tabularnewline
\hline
\end{tabular}\end{gather*}Here $R$ is the root system, $G_{\RR}$ is the simply connected group (which has all $\tilde{V}^{\l}$ as representations), and $d=\dim_{\CC}D=\dim_{\CC}\mathrm{X}$. (See \cite{LZ} for more details.)
\begin{Remark} We have omitted some cases to avoid redundancy due to exceptional isomorphisms and conjugate-isomorphisms. The latter are induced by the action of $\tau$ on $\sigma_{{\tt I}}$. For $A_{n}$, $\tau$ exchanges $\sigma_{i}$ and $\sigma_{n+1-i}$ $(\forall\, i)$; for $D_{n}$, $n$ odd, $\tau$ exchanges $\sigma_{n-1}$ and $\sigma_{n}$; for $E_{6}$, $\tau$ exchanges $\sigma_{6}\leftrightarrow\sigma_{1}$ and $\sigma_{5}\leftrightarrow\sigma_{3}$; and in all other cases the action is trivial.

The exceptional isomorphisms are $\br{III}_{1}\cong\br{I}_{1,1}\cong\br{IV}_{1}\cong\br{II}_{2}$, $\br{III}_{2}\cong\br{IV}_{3}$, $\br{I}_{2,2}\cong\br{IV}_{4}$, $\br{II}_{3}\cong\br{I}_{1,3}$, $\br{II}_{4}\cong\br{IV}_{6}$ (triality for $D_{4}$), and $\br{IV}_{2}\cong\br{III}_{1}\times\br{III}_{1}$. For this reason we consider only $A_{n\geq2}$, $B_{n\geq3}$, $C_{n\geq1}$, and $D_{n\geq4}$.
\end{Remark}
In each case, the real variations arising from $H^{1}$ of the canonical families of abelian varieties over $\mathrm{X}$ are (half-twists of\footnote{This is only relevant for $A_{n}$, cf.~Section~\ref{S4}.}) the $\tilde{\V}^{\l}$ with $-\tfrac{1}{2} ({\tt E}(\l)+{\tt E}(\tau(\l)) )=1$ (cf.~\eqref{eq 6.3}). All such $\l$ take the form $\w_{i}$, with the possibilities corresponding to so-called ``symplectic nodes'' (cf.~\cite{LZ}). This will be recalled case by case where relevant in Sections~\ref{S3} and~\ref{S4}. The variations of CY type are even simpler to describe: they are (again, up to half-twist for $A_{n}$, $D_{n\text{ odd}}$, $E_{6}$) precisely the $\tilde{\V}^{k\w_{{\tt I}}}$ for $k\geq1$~\cite{FL}.

\subsection{The main calculation} \label{S2.4}

We now apply Proposition~\ref{Kostant} to compute the $\big\{\H^{k}(j)\big\}$ for the Hermitian variations $\tilde{\V}^{\l}$ (over $\mathrm{X}$ or $D$) with ${\tt E}(\l)$ odd. Of course, we can work with (its summands) the \emph{complex} variations $\V_{\CC}^{\l}$, and do the computation at one point $\varphi_{0}\in D$. Write\begin{gather} \label{eq 8.1}
W^0(k,j):=\left\{ w\in W^0(k) \left| \tfrac{1}{2}\left( {\tt{E}}(w\cdot \l)-1\right) = j\right. \right\} .
\end{gather}
\begin{Proposition}\label{main}For the $\V_{\CC}^{\l}$, $\H^{k}(j)|_{\varphi_{0}}\cong\oplus_{w\in W^{0}(k,j)}V_{0}^{w\cdot\l}$.\end{Proposition}
\begin{proof}First, we note that \begin{gather} \label{eq 8.2}
\oplus_j \H^k(j) |_{\varphi_0} \cong H^k \big(\fn,V^{\l}\big)
\end{gather}follows from the identification of the Lie algebra cohomology complex $\big(\wedge^{\b}\fn^{\vee}\otimes V^{\l},d\big)$ with the associated graded $\big(\oplus_{j}{\rm Gr}_{F}^{j}C^{\b},\nab\big)$ under the isomorphism $\fn^{\vee}\cong\Omega_{D,\varphi_{0}}^{1}$. To see this identification locally, extend $v\in\big(V^{\l}\big)^{j,-j-1}$ to a section of $\big(\V^{\l}\big)^{j,-j-1}$ by $\tilde{v}=\rho^{\l}(g).v$ ($g\in \exp(\fn)$), and write (for $X\in\fn$) $(\nabla_{X}\tilde{v})|_{\varphi_{0}}=\lim\limits_{\epsilon\to0}\tfrac{1}{\epsilon}\big(e^{\epsilon d\rho^{\l}(X)}v-v\big)=d\rho^{\l}(X).v$
($=$``$X(v)$'').

Next, we compare Hodge gradings on the two sides of \eqref{eq 8.2}. For $X_i^{*}\in\fn^{\vee}$ ($i=1,\ldots,k$) and $v\in {\rm Gr}_{F}^{j-k}V^{\l}=\big(V^{\l}\big)^{j-k,-j+k-1}$, we have ${\tt E}(X_i^{*})=2X_i^{*}$ and ${\tt E}(v)=(2j-2k+1)v$ $\implies$ ${\tt E}\big( \big(X_1^{*}\wedge\cdots\wedge X_k^{*}\big) \otimes v \big)=(2j+1)\big(X_1^{*}\wedge\cdots\wedge X_k^{*} \big) \otimes v$. So the image of $\H^{k}(j)|_{\varphi_{0}}$ in $H^{k}\big(\fn,V^{\l}\big)=\oplus_{w\in W^{0}(k)}V_{0}^{w\cdot\l}$ consists of all the weight spaces with weights $\xi$ for which ${\tt E}(\xi)=2j+1$. Since $\fg^{0}$ commutes with ${\tt E}$, these are just the $V_{0}^{w\cdot\l}$ for which ${\tt E}(w\cdot\l)=2j+1$. The result follows.
\end{proof}

Obviously $W^{0}(0)=\{{\rm id}\}$, and always ${\tt E}(\l)\leq-1$; so immediately we have $\H^{0}(j)=\{0\}$ for $j\geq0$. We claim that $W^{0}(1)=\{{\tt s}\}$, where ${\tt s}$ is the reflection in the distinguished simple root $\sigma_{{\tt I}}$. Indeed, elements of $W^{0}(1)$ have length one, so must be the reflection in one simple root. But the reflection $w_{i}$ in $\sigma_{i\neq{\tt I}}$ doesn't satisfy $w_{i}(\Delta^{+})\supseteq\Delta_{0}^{+}$ (as $\sigma_{i}\notin w_{i}(\Delta^{+})$), whereas ${\tt s}(\Delta_{0}^{+})\subseteq\Delta^{+}$ since ${\tt s}(\sigma_{{\tt I}})=-\sigma_{{\tt I}}$ and $|{\tt s}(\Delta^{+})\cap\Delta^{-}|=1$. So we have the basic
invariant\begin{gather*} 
\mu(\l):=\tfrac{1}{2} ( {\tt{E}} ( {\tt{s}} \cdot \l)-1 ),
\end{gather*}which depends upon the choice of $\sigma_{{\tt I}}$ (hence $D$ and~$\mathrm{X}$). Referring to Section~\ref{S2.2} for the definition of $\tau$, Propositions~\ref{Nori} and~\ref{main} yield at
once the
\begin{Theorem}\label{thm0}If $\mu(\l),\mu(\tau(\l))<0$, then $\operatorname{NF}_{\CT}(\jmath^{*}\V)=\{0\}$ for any \'etale neighborhood \linebreak \smash{$\jmath\colon \CT\to\mathrm{X}$} and any variation $\V$ on $\mathrm{X}$ with the Hermitian variation $\tilde{\V}_{\RR}^{\l}$ as underlying $\RR$-VHS. $($In particular, we have $\overline{\operatorname{AJ}}(Z)=0$ if $\V$, $Z$ are as in~\eqref{eq 3.1}, with $\CS=\mathrm{X}$.$)$
\end{Theorem}

One can of course replace $\mathrm{X}$ by a Zariski open subset in these statements.
\begin{Remark}If $\mu(\l)=\mu(\tau(\l))=0$, then (arguing as in~\cite{Fa}) one can show that M. Saito's canonical mixed Hodge structure \cite{Sa} on $H^{1}(\mathrm{X},\VV)$ is pure of type $(0,0)$. (Naturally, $H^{1}(\mathrm{X},\VV)$ could still be $\{0\}$.)
\end{Remark}

Suppose that for $\mathcal{V}$ as above we have $\mu(\l)\geq0$, so that $\mathcal{E}_{1}^{1-\mu(\l),\mu(\l)}=\H_{\nab}^{1}(\mu(\l))\neq\{0\}$ (cf.~\eqref{eq 2.2}). To conclude the stronger result that $\mathcal{E}_{\infty}^{1-\mu(\l),\mu(\l)}$ hence $\H_{\nabla}^{1}\big(F^{0}C^{\b}\big)$ is nonzero, we would need to compute the nonlinear $d_{\ell}\colon \mathcal{E}_{\ell}^{1-\mu(\l),\mu(\l)}\to\mathcal{E}_{\ell}^{2-\mu(\l)-\ell,\mu(l)+\ell}$. The present methods are only of use here if they show that all $\H_{\nab}^{2}(j)=\{0\}$ for $j>\mu(\l)$. But this is asking too much: even in the key case $\{D=\br{III}_{3},\, \l=\w_{3},\, \mu(\l)=0\}$ where we have a nonzero \emph{geometric} normal function coming from the Ceresa cycle, we compute that $\H^{2}(1)\neq\{0\}$. So with regard to predicting normal functions it is probably more useful to stick to $\H^{1}(j)$ and make the provisional
\begin{Definition} \label{DefINF} The complex (resp.\ real) variation $\V^{\l}$ (resp.~$\tilde{\V}_{\RR}^{\l}$) has an \emph{infinitesimal normal function} if $\mu(\l)\geq0$.
\end{Definition}

\section{Analysis of (mostly) tube domain cases} \label{S3}

In this section we study the real variations $\tilde{\V}_{\RR}^{\l}$ arising from dominant integral $\l$ with ${\tt E}(\l)$ an odd integer. Though we (more or less) carry this out for all the domains in the table of Section~\ref{S2.3}, the results are of interest mainly when $D$ is of tube type ($\br{I}_{p,p}$, $\br{II}_{2m}$, $\br{III}_{n}$, $\br{IV}_{m}$, or~$\br{EVIII}$). As above, we have ${\tt E}(\sigma_{{\tt I}})=-2$, ${\tt E}(\sigma_{j\neq{\tt I}})=0$, and write $\l=\sum\limits_{i=1}^{n}m_{i}\w_{i}$; note that ${\tt s}$ sends $\sigma_{{\tt I}}\mapsto-\sigma_{{\tt I}}$ and fixes all $\w_{j\neq{\tt I}}$. To streamline the discussion of examples, we make the
\begin{Definition}A family $\CA\overset{\pi}{\to}\CS$ of abelian varieties \emph{admits no reduced normal functions} if $\operatorname{NF}_{\CS}(\V)=\{0\}$ for every $r\in\mathbb{N}$ and irreducible VHS $\V\subset R^{2r-1}\pi_{*}\QQ_{\CA}(r)\otimes\mathcal{O}_{\CS}$ of level $>1$.
\end{Definition}
In all the cases below where we obtain such an assertion, one knows that all level-one real sub-VHS are in fact defined over $\QQ$.

\subsection[Case $\br{III}_{n}$]{Case $\boldsymbol{\br{III}_{n}}$} \label{S3.1}

This is the case analyzed by Nori \cite{No} ($\H^{1}(j)$) and Fakhruddin~\cite{Fa} ($\H^{k}(j)$ for $k>1$). We have
\begin{gather*}
\Sigma=\{e_{1}-e_{2},e_{2}-e_{3},\ldots,e_{n-1}-e_{n},2e_{n}\},\\
\Omega=\{e_{1},e_{1}+e_{2},e_{1}+e_{2}+e_{3},\ldots,e_{1}+\cdots+e_{n}\},
\end{gather*}
and so ${\tt I}=n$ $\implies$ ${\tt E}(\l)=-\sum\limits_{i=1}^{n}im_{i}$; the only level-one variation is thus $\V_{\RR}^{\w_{1}}$ (which arises from~$H^{1}$ of the universal abelian variety). There are no complex or quaternionic irreps, and (taking~${\tt E}(\l)$ odd) we have ${\tt s}(\w_{n})=2\w_{n-1}-\w_{n}$ $\implies$
\begin{gather*}
\mu(\l)=-\tfrac{1}{2}\sum_{i<n}im_{i}+\big(1-\tfrac{n}{2}\big)m_{n}+\tfrac{1}{2}.
\end{gather*}
So for $\ell\big(V^{\l}\big)=-{\tt E}(\l)>1$ odd, $\mu(\l)\geq0$ for only
\begin{gather*}
n=1\quad \text{and} \quad \l=a\w_{1}, \quad a>1\text{ odd},\\
n=2 \quad \text{and} \quad \l=\w_{1}+a\w_{2}, \quad a\geq1,\\
n=3\quad \text{and} \quad \l=\w_{3};
\end{gather*}
in particular, no \'etale pullback of the universal abelian variety $\CA^{\br{III}_{g}}\to\A_{g}(\Gamma)=\Gamma\backslash\br{III}_{g}$ admits reduced normal functions for $g\geq4$. These are the results of Nori~\cite{No}.

\subsection[Case $\br{EVIII}$]{Case $\boldsymbol{\br{EVIII}}$} \label{S3.2}

We have ${\tt I}=7$ and
\begin{gather*}
\Sigma=\big\{ \tfrac{1}{2}(e_{8}-e_{7}-\cdots-e_{2}+e_{1}),e_{2}+e_{1},e_{2}-e_{1},e_{3}-e_{2},\ldots,e_{6}-e_{5}\big\} ;
\end{gather*}
$\Omega$ is written in terms of the $\{\sigma_{i}\}$ in \cite[p.~688]{Kn}. This yields
\begin{gather*}
{\tt E}(\l)=-2m_{1}-3m_{2}-4m_{3}-6m_{4}-5m_{5}-4m_{6}-3m_{7},
\end{gather*}
${\tt E}({\tt s}(\w_{7}))=-1$, and
\begin{gather*}
\mu(\l)=-m_{1}-\tfrac{3}{2}m_{2}-2m_{3}-3m_{4}-\tfrac{5}{2}m_{5}-2m_{6}-\tfrac{1}{2}m_{7}+\tfrac{1}{2}.
\end{gather*}
There are no complex or quaternionic irreps, and $\mu(\l)\geq0$ (with ${\tt E}(\l)$ odd) for
\begin{gather*}
\l=\w_{7}.
\end{gather*}

\subsection[Case $\br{EIII}$]{Case $\boldsymbol{\br{EIII}}$} \label{S3.3}

Here ${\tt I}=1$, $\Omega$ is expressed in terms of
\begin{gather*}
\Sigma=\big\{ \tfrac{1}{2}(e_{8}-e_{7}-\cdots-e_{2}+e_{1}),e_{1}+e_{2},e_{2}-e_{1},\ldots e_{5}-e_{4}\big\}
\end{gather*}
in \cite[p.~687]{Kn}, and $\tau=(16)(35)$ on both the $\{\sigma_{i}\}$ and $\{\w_{i}\}$. Since
\begin{gather*}
{\tt E}(\l)=-\tfrac{2}{3}(4m_{1}+5m_{3}+4m_{5}+2m_{6})-2m_{2}-4m_{4},
\end{gather*}
we find that ${\tt E}(\l)\in\ZZ$ $\iff$ $3|(4m_{1}+5m_{3}+4m_{5}+2m_{6})$ $\iff$ ${\tt E}(\l)\in2\ZZ$. So there are no odd-level VHS without half-twists, and further analysis is deferred to Section~\ref{S4}.

We can now turn to the richer remaining classical cases.

\subsection[Case $\br{II}_{n}$]{Case $\boldsymbol{\br{II}_{n}}$} \label{S3.4}

With ${\tt I}=n$ (note that $n\geq4$),
\begin{gather*}
\Sigma=\{e_{1}-e_{2},\ldots,e_{n-1}-e_{n},e_{n-1}+e_{n}\},\\
\Omega=\big\{e_{1},e_{1}+e_{2},\ldots,e_{1}+\cdots+e_{n-2},\tfrac{1}{2}(e_{1}+\cdots+e_{n-1}-e_{n}),\tfrac{1}{2}(e_{1}+\cdots e_{n})\big\}
\end{gather*}
we obtain
\begin{gather*}
{\tt E}(\l)=-\sum_{i=1}^{n-2}im_{i}-\tfrac{n-2}{2}m_{n-1}-\tfrac{n}{2}m_{n}.
\end{gather*}
Note that $\tau$ swaps $\sigma_{n-1}\leftrightarrow\sigma_{n}$ and $\w_{n-1}\leftrightarrow\w_{n}$ if $n$ is odd, and is trivial for~$n$ even. The only complex irreps are therefore the $V^{\l}$, with $n$ odd and $m_{n-1}\neq m_{n}$. The quaternionic ones are the~$V^{\l}$ with
\begin{gather*}
\sum_{i=1}^{n-2}im_{i}+m_{n-1} \quad \text{odd (}n\text{ even)},\qquad \text{or}\qquad \begin{cases}
\displaystyle \sum_{i=1}^{n-2}im_{i}\text{ odd}\\
m_{n-1}=m_{n}
\end{cases} \quad (n\text{ odd}).
\end{gather*}
We have $\ell\big(\tilde{V}^{\l}\big)=1$ for $\l=\w_{1},\w_{3}$ ($n=4$) and $\l=\w_{1}$ ($n>4$); these $V^{\l}$ are all quaternionic (i.e., $\tilde{V}^{\l}=\big(V^{\l}\big)^{\oplus2}$ over~$\CC$).
\begin{Definition} \label{DefQAV} By a \emph{universal quaternionic abelian variety}, we shall mean any family $\CA^{\br{II}_{n}}\to\Gamma\backslash\br{II}_{n}$ of abelian $2n$-folds whose $H^{1}$ recovers $\tilde{\V}_{\RR}^{\w_{1}}$.
\end{Definition}
Such families admit an embedding of a definite rational quaternion algebra $\mathcal{Q}$ into $\operatorname{End}(\CA)_{\QQ}$. (The Mumford--Tate group $G$ is a $\QQ$-form of $G_{\RR}$ and so, by considering its fixed 2-tensors, $\mathcal{Q}$ is a $\QQ$-form of~$\mathbb{H}$. See~\cite{vGV} for more details on quaternionic abelian varieties.) There are natural embeddings $\br{II}_{n}\hookrightarrow\br{III}_{2n}$ which yield countably many ``quaternionic subfamilies'' of $\CA^{\br{III}_{2n}}$.

Now from ${\tt s}(\w_{n})=\w_{n-2}-\w_{n}$ we find
\begin{gather*}
\mu(\l)=-\tfrac{1}{2}\sum_{i=1}^{n-2}im_{i}-\tfrac{n-2}{4}m_{n-1}-\tfrac{n-4}{4}m_{n}+\tfrac{1}{2}.
\end{gather*}
Imposing $-{\tt E}(\l)$ odd $>1$, $n\geq4$, and $\mu(\l)\geq0$ yields
\begin{gather*}
\l=\w_{1}+a\w_{4},\qquad \w_{3}+a\w_{4}, \qquad n=4,\qquad a>0,
\end{gather*}
which are quaternionic, and
\begin{gather*}
\l=\w_{6}, \qquad n=6,
\end{gather*}
which is real.
\begin{Theorem}\label{thm1}No \'etale pullback of a universal quaternionic abelian variety of $($relative$)$ dimension $2n$ admits reduced normal functions outside the case $n=4$.\end{Theorem}
\begin{proof}It remains to deal with $n=6$. The point is that only $\V^{\w_{5}+\w_{6}}$ and $\V^{2\w_{5}}$, $\V^{2\w_{6}}$ occur in $H_{\rm rel}^{*}\big(\CA^{\br{II}_{6}}\big)\cong\bigwedge^{*}\big(\V^{\w_{1}}\oplus\V^{\w_{1}}\big)$; the half-spin variations\footnote{In all the $D_{n}$ cases, $V^{\w_{n-1}}$ and $V^{\w_{n}}$ identify with ${\rm spin}^{-}$ and ${\rm spin}^{+}$ (with order depending on parity of $n$). These are representations of $\operatorname{Spin}^{*}(12)$ but not of ${\rm SO}^{*}(12)$, the Mumford--Tate group of $\CA^{\br{II}_{6}}$.} $\V^{\w_{5}}$, $\V^{\w_{6}}$ do not. (See~\cite{FH} for the equivalent fact on representations of~${\rm SO}(12)$.)
\end{proof}

What is special about the quaternionic $8$-folds that might yield $\overline{\operatorname{AJ}}$-nontrivial elements of the Griffiths group? According to~\cite{vGV}, there exist families $\CA^{\br{II}_{4}}$ whose general member arises as a quaternionic Prym variety, associated to a certain $8:1$ unramified cover of a general genus $3$ curve (plausible as $\dim\M_{3}=\dim\br{II}_{4}=6$). Since $\bigwedge^{3}V^{\w_{1}}=V^{\w_{3}+\w_{4}}$, the (non-CY) variation~$\tilde{\V}_{\RR}^{\w_{3}+\w_{4}}$ occurs in $H^{3}$ of $\CA^{\br{II}_{4}}$, and it seems reasonable to expect that one can construct $\overline{\operatorname{AJ}}$-nontrivial 1-cycles from the Abel--Prym image of the genus $17$ cover curves. For $n=6$, it seems to be an open problem to give a simple motivic construction of $\V_{\RR}^{\w_{6}}$, so we cannot speculate about cycles in this case.

\subsection[Case $\br{I}_{p,n-p+1}$]{Case $\boldsymbol{\br{I}_{p,n-p+1}}$} \label{S3.5}

We have $n\geq2$, $1\leq p\leq\big\lfloor \tfrac{n+1}{2}\big\rfloor $, ${\tt I}=p$,
\begin{gather*}
\Sigma=\{e_{1}-e_{2},e_{2}-e_{3},\ldots,e_{n}-e_{n+1}\}, \\
\Omega=\big\{\w_{i}=e_{1}+\cdots+e_{i}-\tfrac{i}{n+1}(e_{1}+\cdots+e_{n})\big\}_{i=1,\ldots,n},
\end{gather*}
with $\tau=\prod\limits_{j=1}^{\lfloor \frac{n}{2}\rfloor }(j,n+1-j)$ on the $\{\sigma_{i}\}$ and $\{\w_{i}\}$, and ${\tt s}(\w_{p})=\w_{p-1}-\w_{p}+\w_{p+1}$ (or $-\w_{1}+\w_{2}$ if $p=1$). This yields
\begin{gather}
 {\tt E}(\l)=\tfrac{2}{n+1}\left(\sum_{i=1}^{p}(n+1-p)im_{i}+\sum_{i=p+1}^{n}p(n+1-i)m_{i}\right), \nonumber\\
{\tt E}({\tt s}\cdot\l)=2+2m_{p}+{\tt E}(\l),\qquad \mu(\l)=\tfrac{1}{2}+m_{p}+\tfrac{1}{2}{\tt E}(\l)\label{eq 15.1}
\end{gather}and for $n$ odd there are quite a few (mostly complex) representations with ${\tt E}(\l)$ an odd integer and $\mu(\l)\geq0$. For example, $p=2$, $n=7$ produces $\tilde{V}^{\w_{6}}=V^{\w_{2}}\oplus V^{\w_{6}}$ and $\tilde{V}^{2\w_{7}}=V^{2\w_{7}}\oplus V^{2\w_{1}}$, which both have level $3$ (and $\mu(2\w_{7})=\mu(\w_{6})=\mu(\w_{2})=0$). We have not carried out an exhaustive search, because the resulting variations seem obscure and the situation drastically improves with the introduction of half-twists (Section~\ref{S4}).

So restricting to $p=\tfrac{n+1}{2}$ ($n$ odd), i.e., the case $\br{I}_{p,p}$ ($p\geq2$), we have
\begin{gather*}
{\tt{E}}(\l) = -\sum_{i=1}^p im_i - \sum_{i=p+1}^{2p-1} (2p-i) m_i,\nonumber \\
\mu(\l) = -\tfrac{1}{2}\sum_{i<p} im_i - \big(\tfrac{p}{2} - 1\big)m_p - \sum_{i>p} \big(p-\tfrac{i}{2}\big)m_i + \tfrac{1}{2} .
\end{gather*}We remark that all representations are either real ($m_{j}=m_{n+1-j}$ $\forall\, j$) or complex, and that the sole level $1$ variation is $\tilde{\V}^{\w_{1}}=\V^{\w_{1}}\oplus\V^{\w_{n}}$.
\begin{Definition} \label{DefWAV} A \emph{universal Weil abelian variety} is a family $\CA^{\br{I}_{p,p}}\to\Gamma\backslash\br{I}_{p,p}$ of abelian $2p$-folds (of dimension $p^{2}$) whose $H^{1}$ recovers $\tilde{\V}_{\RR}^{\w_{1}}$.
\end{Definition}
Such families admit an embedding of an imaginary quadratic field into $End(\CA)_{\QQ}$, and produce countably many subfamilies of~$\CA^{\br{III}_{2p}}$.
\begin{Theorem} \label{thm2}No \'etale pullback of a universal Weil abelian variety of $($relative$)$ dimension~$2p$ admits reduced normal functions outside the cases $p=2,3$.\end{Theorem}
\begin{proof} For $p>3$, there are no solutions to $\mu(\l)\geq0$ (with the usual constraints).
\end{proof}

Indeed, the only solutions are
\begin{gather*}
\l=\w_{1}+a\w_{2},\qquad \w_{3}+a\w_{2}, \qquad p=2,
\end{gather*}
which are conjugate complex, and
\begin{gather*}
\l=\w_{3}, \qquad p=3,
\end{gather*}
which is real (with $\V^{\w_{3}}$ of CY type). One easily shows that $\tilde{\V}^{\w_{1}+\w_{2}}\subset H_{\rm rel}^{3}\big(\CA^{\br{I}_{2,2}}\big)$ and $\tilde{\V}^{\w_{3}}\subset H_{\rm rel}^{3}\big(\CA^{\br{I}_{3,3}}\big)$,\footnote{Since $H_{\rm rel}^3\big(\CA^{\br{I}_{p,p}}\big)\cong \bigwedge^3 \tilde{\V}^{\w_1} = \bigwedge^3 (\V^{\w_1}\oplus\V^{\w_{2p-1}})$, this follows from the fact that $\bigwedge^k \V^{\w_1}= \V^{\w_k}$ (and $\bigwedge^k \V^{w_{2p-1}}= \V^{w_{2p-k}}$) for $k\leq 2p-1$ and $\V^{\w_1 + \w_2}\subset \V^{\w_1}\otimes\V^{\w_2}$. See \cite[Section~15]{FH}.}
and these are the two cases where there exist families $\CA^{\br{I}_{p,p}}$ with generic member arising from a generalized (3~to~1) Prym construction. More precisely, let $\mathcal{B}_{p,p}$ denote the moduli space of abelian $2p$-folds $A$ with: (a)~a degree-$3$ automorphism $\alpha$, having $\zeta_{3}$- and $\bar{\zeta}_{3}$-eigenspaces in $T_{0}A$ both of dimension~$p$; and (b)~a~polarization $\Sigma$ with $\ker(\Sigma)\subset A[3]$ of order~$3^{2p}$. Then denoting by $\mathcal{R}_{g}^{3}$ the moduli space of unramified $3:1$ curve covers $\tilde{C}\to C$ with $g_{C}=g$, the Prym map $\mathcal{P}\colon \mathcal{R}_{g}^{3}\to\mathcal{B}_{g-1,g-1}$ defined by\footnote{Note that via $\operatorname{pr} :=2{\rm id}-\alpha_{*}-\alpha_{*}^{2}\colon J(\tilde{C})\twoheadrightarrow A$, we may view $A$ as the quotient $J(\tilde{C})/J(C)$.} $A:=\ker\big(J\big(\tilde{C}\big)\twoheadrightarrow J(C)\big)^{\circ}$ is dominant for $g=3,4$~\cite{F}. This was used by Schoen~\cite{Sc} to prove the Hodge conjecture for the specific Weil~4- and 6-folds parametrized by $\mathcal{B}_{2,2}$ and $\mathcal{B}_{3,3}$.

Let $\tilde{C}^{+}$ [resp.\ $\tilde{C}^{-}$] denote the Abel image of $\tilde{C}$ in $J\big(\tilde{C}\big)$ [resp.\ its image under $-{\rm id}$], defined up to algebraic equivalence, and write $\mathcal{O}\subseteq\mathcal{O}_{\mathbb{Q}(\sqrt{-3})}$ for an order. We now demonstrate that in the $g_{C}=4$ ($g_{\tilde{C}}=10$) case, the \emph{Prym--Ceresa cycle}
\begin{gather*}
Z:=\operatorname{pr}_{*}\big(\tilde{C}^{+}-\tilde{C}^{-}\big)\in\operatorname{Griff}^{5}(A)_{\mathbb{Q}}
\end{gather*}
generates a nontrivial reduced normal function on $\mathcal{R}_{4}^{3}(\mathcal{O})$ (over $\mathcal{B}_{3,3}(\mathcal{O})$).\footnote{Moduli spaces with $\mathcal{O}$-level structure.} (This is in contrast to the more standard $2:1$ Prym setting, where such cycles are algebraically equivalent to zero.) This is accomplished by degenerating $C$ to a union $D\cup E$ ($D$ of genus $3$, $E$ elliptic, meeting at a single node), so that the cover $\tilde{C}$ becomes $(D_{1}\sqcup D_{2}\sqcup D_{3})\cup\tilde{E}$ (each $D_{i}\cong D$ meeting $\tilde{E}\simeq E$ once). The cyclic automorphism $\alpha$ permutes the $\{D_{i}\}$, and the specialized Prym variety $A$ is the quotient of $J(E)\oplus\bigoplus_{i=1}^{3}J(D_{i})$ by the invariant part $J(E)\oplus J(D)_{\Delta}$ ($\Delta=$ diagonal). We have $H^{1}(A)\cong V^{\w_{1}}\oplus V^{\w_{5}}$, and $H^{3}(A)^{\alpha}\cong V^{\w_{3}}\oplus V^{\w_{3}}$. Writing $\big\{\w_{i}^{j}\big\}_{j=1,2,3}\subset\Omega^{1}(J(D_{i}))$ for bases, the invariant form
\begin{gather*}
\w:=\big(\w_{1}^{1}+\zeta_{3}\w_{2}^{1}+\bar{\zeta}_{3}\w_{3}^{1}\big)\wedge\big(\w_{1}^{2} +\zeta_{3}\w_{2}^{2}+\bar{\zeta}_{2}\w_{3}^{2}\big)\wedge\big(\w_{1}^{3}+\zeta_{3}\w_{2}^{3}+\bar{\zeta}_{3}\w_{3}^{3}\big)
\end{gather*}
is the pullback of a class from $H^{3,0}(A)^{\alpha}$. The projection of the Ceresa cycle $D_{1}^{+}-D_{1}^{-}=\partial\Gamma$ (for $D_{1}$) to~$A$ has $\int_{\Gamma}\omega=\int_{\Gamma}\w_{1}^{1}\wedge\w_{1}^{2}\wedge\w_{1}^{3}$ generically nonzero (otherwise $\overline{\operatorname{AJ}}$ would not detect the usual generic Ceresa cycle). But the degeneration of the Prym--Ceresa cycle is the sum of the projections of the $D_{i}^{+}-D_{i}^{-}$ to $A$, and (by invariance under $\alpha$) each has the same nonzero image under the composition of $\overline{\operatorname{AJ}}^{5}$ with the projection $J\big(H^{3}(A)^{\vee}\big)\twoheadrightarrow J\big(\big(H^{3}(A)^{\alpha}\big)^{\vee}\big)$. So we have the desired geometric realization of $(\br{I}_{3,3},\omega_{3})$.

Applying the (algebraic) Lefschetz operator to $Z$, we obtain cycles (over $\mathcal{R}_{4}^{3}(\mathcal{O})$) $Z^{k}\in\operatorname{Griff}^{k}(A)_{\mathbb{Q}}$ for $k=2,3,4,5$ with nontorsion $\overline{\operatorname{AJ}}^{k}$ image. By \cite{Ra} (see Section~\ref{S1} above), $\mathcal{R}_{4}^{3}(\mathcal{O})$\footnote{More precisely, we get a normal function over the minimal quotient of $\mathcal{R}_{4}^{3}(\mathcal{O})$ over $\mathcal{B}_{3,3}(\mathcal{O})$ which supports $Z$ (modulo algebraic equivalence).} cannot be an $\br{I}_{3,3}(\Gamma)$ ($\Gamma$ arithmetic), and so has a nonempty branch locus over $\mathcal{B}_{3,3}(\mathcal{O})$. Moreover, the push-forward of $\overline{\operatorname{AJ}}^{k}(Z^{k})$ to $\mathcal{B}_{3,3}(\mathcal{O})$ must be zero, so that the Galois group of the cover acts nontrivially. Now the last paragraph of \cite{No2} applies, with the obvious substitutions: we push forward $Z^{k}$ along the isogenies of the universal Weil $6$-fold over the automorphisms of~$\br{I}_{3,3}$ induced by~$G(\mathbb{Q})$ (regarding this as taking place in the inverse limit over level structures~$\mathcal{O}$). Since the action of $G(\mathbb{Q})$ produces infinitely many translates of the original base locus, the absolute Galois group of $\underset{\mathcal{O}}{\underrightarrow{\lim}}\,\mathbb{C}(\mathcal{B}_{3,3}(\mathcal{O}))$ acts through infinitely many distinct representations upon the push-forwards of~$Z^{k}$. So we have the
\begin{Theorem}\label{thm Griff Weil} Let $A$ be a very general Weil abelian $6$-fold in the family parametrized by $\mathcal{B}_{3,3}$. Then the groups $\operatorname{Griff}^{k}(A)_{\mathbb{Q}}$, and their $\overline{\operatorname{AJ}}^{k}$-images, are countably infinite-dimensional for $k=2,3,4,5$.
\end{Theorem}

\subsection[Cases $\br{IV}_{2n-2}$ and $\br{IV}_{2n-1}$]{Cases $\boldsymbol{\br{IV}_{2n-2}}$ and $\boldsymbol{\br{IV}_{2n-1}}$} \label{S3.6}

Beginning with $\br{IV}_{2n-2}$ (the $D_{n}$ case), we note that $n\geq4$, ${\tt I}=1$, $\Sigma$ and $\Omega$ are as in Section~\ref{S3.4}, and
\begin{gather*}
{\tt E}(\l) (={\tt E}(\tau(\l)) )=-2(m_{1}+\cdots+m_{n-2})-(m_{n-1}+m_{n}).
\end{gather*}
$V^{\l}$ is complex iff $n$ is odd and $m_{n-1}\neq m_{n}$, and quaternionic iff $4|n$ and $m_{n-1}+m_{n}$ is odd. So the level-one real Hermitian VHS are ($\otimes\CC$)
\begin{alignat*}{4}
& (n\text{ odd}) \qquad && \tilde{\V}^{\w_{n}}=\V^{\w_{n}}\oplus\V^{\w_{n-1}},\qquad && &\\
& (4|n) \qquad && \tilde{\V}^{\w_{n}}=\big(\V^{\w_{_{n}}}\big)^{\oplus2}, \qquad && \tilde{\V}^{\w_{n-1}}=\big(\V^{\w_{n-1}}\big)^{\oplus2}, &\\
& (4|n+2) \qquad & & \tilde{\V}^{\w_{n}}=\V^{\w_{n}}, \qquad && \tilde{\V}^{\w_{n-1}}=\V^{\w_{n-1}}.&
\end{alignat*}
From ${\tt s}(\w_{1})=\w_{2}-\w_{1}$, we define
\begin{gather*}
\mu(\l)=\tfrac{1}{2}-m_{2}-\cdots-m_{n-2}-\tfrac{1}{2}m_{n-1}-\tfrac{1}{2}m_{n},
\end{gather*}
so that $-{\tt E}(\l)>1$ odd and $\mu(\l)\geq0$ $\implies$
\begin{gather*}
\l=a\w_{1}+\w_{n-1}\qquad \text{or}\qquad a\w_{1}+\w_{n}, \qquad a>0,
\end{gather*}
which works for arbitrary $n$!

A similar (but simpler) story unfolds in the $B_{n}$ case $\br{IV}_{2n-1}$, where (with $n\geq3$, ${\tt I}=1$)
\begin{gather*}
\Sigma=\{e_{1}-e_{2},\ldots,e_{n-1}-e_{n},e_{n}\}, \\
\Omega=\big\{e_{1},e_{1}+e_{2},\ldots,e_{1}+\cdots+e_{n-1},\tfrac{1}{2}(e_{1}+\cdots+e_{n})\big\},\\
{\tt E}(\l)=-2(m_{1}+\cdots+m_{n-1})-m_{n},
\end{gather*}
and (with no complex irreps, as $\tau={\rm id}$) the only quaternionic irreps are the $V^{(2k+1)\w_{n}}$ for $\big\lfloor \tfrac{n-1}{2}\big\rfloor $ odd. The level one Hermitian $\RR$-VHS are thus ($\otimes\CC$) the
\begin{gather*}
\tilde{\V}^{\w_{n}}=\big({\V}^{\w_{n}}\big)^{\oplus2}, \quad n=3,4;7,8;\ldots \qquad \text{or}\qquad \V^{\w_{n}}\quad \text{otherwise}.
\end{gather*}
We have ${\tt s}(\w_{1})=\w_{2}-\w_{1}$, hence
\begin{gather*}
\mu(\l)=\tfrac{1}{2}-m_{2}-\cdots-m_{n-1}-\tfrac{1}{2}m_{n},
\end{gather*}
so that the variations associated to
\begin{gather*}
\l=a\w_{1}+\w_{n}, \qquad a>0,
\end{gather*}
are (for any $n$) the ones with infinitesimal reduced normal functions.
\begin{Definition}A \emph{universal spin abelian variety} is an abelian family $\CA\to\Gamma\backslash\br{IV}_{m}$ with $H^{1}$ (over $\RR$) a number of copies of $\tilde{\V}_{\RR}^{\w_{n-1}}$ and $\tilde{\V}_{\RR}^{\w_{n}}$ ($m=2n-2$ even) resp. $\tilde{\V}_{\RR}^{\w_{n}}$ ($m=2n-1$ odd).
\end{Definition}
The minimum possible (relative) dimension of $\CA$ is clearly $2^{n}$ ($m$ and $\big\lfloor \tfrac{n-1}{2}\big\rfloor $ odd), $2^{n-1}$ ($m$~and~$\big\lfloor \tfrac{n+1}{2}\big\rfloor $ odd; $m$ even and $4\nmid n+2$), or~$2^{n-2}$ ($m$~even and~$4|n+2$). However, the main natural source of spin abelian varieties is the Kuga--Satake construction (cf.~\cite{vG2}), which produces varieties of much higher dimension: for $m\leq19$, one has families of $K3$ surfaces $\CX$ with $H_{\rm tr}^{2}(\CX)\cong\V^{\w_{1}}$, and Clifford algebras produce an embedding $H_{\rm tr}^{2}(\CX)\hookrightarrow H^{1}(\CA)^{\otimes2}$ with
\begin{gather*}
H^{1}(\CA)\cong \begin{cases}
\big(\V^{\w_{n-1}}\big)^{\oplus2^{n-1}}\oplus\big(\V^{\w_{n}}\big)^{\oplus2^{n-1}}, & m\text{ even},\\
\big(\V^{\w_{n}}\big)^{\oplus2^{n}}, & m\text{ odd}.
\end{cases}
\end{gather*}
So the situation is in marked contrast to those encountered above:\footnote{We haven't ruled out that for some ``minimal'' spin abelian varieties one might have a vanishing result analogous to the above ones, except for $\br{IV}_{m}$ with $m=7,8,9;13,15$ etc.}
\begin{Proposition}The relative cohomology of any Kuga--Satake family of spin abelian varieties over $\br{IV}_{m}$ $($any $m\geq7)$ admits infinitesimal $($reduced$)$ normal functions.\end{Proposition}

\begin{proof}We only need to show that (say) $\V^{\w_{1}+\w_{n}}$ occurs in $H^{3}(\CA)=\bigwedge^{3}H^{1}(\CA)$. This is done by considering the decomposition of $\big(V^{\w_{n}}\big)^{\otimes2}$ or $V^{\w_{n}}\otimes V^{\w_{n-1}}$; e.g., for $m$ odd,
\begin{gather*}
\big(V^{\w_{n}}\big)^{\otimes2}\cong V^{2\w_{n}}\oplus V^{\w_{n-1}}\oplus\cdots\oplus V^{\w_{1}}\oplus\mathbf{1}
\ \implies\ \big(V^{\w_{n}}\big)^{\otimes3}\supseteq V^{\w_{1}}\otimes V^{\w_{n}}\supseteq V^{\w_{1}+\w_{n}}.
\end{gather*}
(There are enough copies of $\V^{\w_{n}}$ that we can consider tensor rather than wedge powers.)
\end{proof}

Finally, we remark that the triality isomorphism for $D_{4}$ exhibits the universal quaternionic abelian $8$-folds as ``minimal'' spin $8$-folds, by equating $\tilde{\V}_{\RR}^{\w_{4}}\to\br{IV}_{6}$ and $\tilde{\V}_{\RR}^{\w_{1}}\to\br{II}_{4}$. Specialization under $\br{IV}_{5} \hookrightarrow\br{IV}_{6}$ shrinks the Mumford--Tate group to $\operatorname{Spin}(2,5)$ and gives geometric realizations of $\tilde{\V}_{\RR}^{\w_{3}} \to\br{IV}_{5}$ (as noticed by~\cite{vGV}). All these cases admit infinitesimal normal functions.

\section{Half-twists and non-tube cases} \label{S4}

In this section we consider a slight generalization of the homogeneous variations described in Section~\ref{S2.3}, by enlarging our simple $G_{\CC}$ to $\tilde{G}_{\CC}=\mathbb{G}_{m}\cdot G_{\CC}$ (or $G_{\RR}$ to~$U(1)\cdot G_{\RR}$) and taking~$\tilde{{\tt E}}$ to have a component in the abelian factor of $\tilde{\fg}$. This allows us to shift the Hodge grading of a complex summand (e.g., to make it integral), which is necessary in order to study the cohomology of abelian varieties of generalized Weil type and to obtain all CY variations. We restrict our investigation to these examples and proceed with a~minimum of formality.

Given an irrep $V^{\l}$ of $\fg$ and ${\tt E}\in\ft$ as before, we take $\tilde{{\tt E}}=({\tt E},1)\in\fg\oplus\CC=\tilde{\fg}$, and define representations of~$\tilde{\fg}$ by
\begin{gather*}
\big( V^{\l}\big\{ \tfrac{a}{2}\big\} \big)^{p,-p-1}:=\big(V^{\l}\big)^{p+\frac{a}{2},-p-\frac{a}{2}-1}.
\end{gather*}
We write\begin{gather} \label{eq 20.1}
\tilde{V}^{\l}\big\{ \tfrac{a}{2}\big\} := V^{\l}\big\{ \tfrac{a}{2}\big\} \oplus V^{\tau(\l)}\big\{\tfrac{-a}{2}\big\} ,
\end{gather}
whether or not $V^{\l}$ or $V^{\l}\big\{ \tfrac{a}{2}\big\} $ is complex; this is a variant of van Geemen's half-twist~\cite{vG1} that preserves the weight. As a real Hodge structure, \eqref{eq 20.1} has level
\begin{gather} \label{eq 20.2}
\ell \big( \tilde{V}^{\l}\big\{ \tfrac{a}{2} \big\} \big) = \max \big\{ {-}{\tt{E}}(\l)+a, -{\tt{E}}(\tau(\l))-a \big\} .
\end{gather}When \eqref{eq 20.2} is odd, define\begin{gather} \label{eq 20.3}
\mu(\l,a):=\tfrac{1}{2}\left\{ {\tt{E}}({\tt{s}}\cdot \l)-a-1\right\} .
\end{gather}
\begin{Proposition}Assume \eqref{eq 20.2} is odd. If either $\mu(\l,a)$ or $\mu(\tau(\l),-a)$ is $\geq0$, then the corresponding weight-$(-1)$ $\mathbb{R}$-VHS $\tilde{\V}_{\RR}^{\l}\big\{ \tfrac{a}{2}\big\} $ over $\mathrm{X}=\Gamma\backslash D$ has a nonzero $\H^{1}(j)$, $j\geq0$. $($As above, we shall say it has an infinitesimal normal function.$)$ If both invariants are negative, any underlying $\QQ$-VHS admits no nontrivial reduced normal function on an \'etale neighborhood of~$\mathrm{X}$.
\end{Proposition}

\begin{proof}Same as that of Proposition~\ref{main} and Theorem \ref{thm0}, but with the condition $\tfrac{1}{2} ({\tt E}(w\cdot \lambda )-1 )=j$ in \eqref{eq 8.1} replaced by $\tfrac{1}{2} ( {\tt E}(w\cdot \lambda)-1 ) = j+\tfrac{a}{2}$ (as $V^{\lambda}$ is replaced by $V^{\lambda}\big\{ \tfrac{a}{2}\big\}$).
\end{proof}

\subsection[Case $\br{I}_{p,n-p+1}$ (bis)]{Case $\boldsymbol{\br{I}_{p,n-p+1}}$ (bis)} \label{S4.1}

As before, we assume $p\leq\big\lfloor \tfrac{n+1}{2}\big\rfloor $. Referring to~\eqref{eq 15.1}, the only choices that make~\eqref{eq 6.3} equal~$1$ are $\l=\w_{1}$ or $\w_{n}$, unless $p=1$ (in which case all $\w_{i}$ work). We shall focus attention on the variations arising in the odd exterior powers of
\begin{gather} \label{eq 21.1}
\H_{p,n-p+1} :=\tilde{\V}^{\w_1}_{\RR}\left\{\frac{2p-n-1}{2(n+1)}\right\}, \qquad \text{i.e.,} \quad a=\tfrac{2p-n-1}{n+1}
\end{gather}
over $\mathrm{X}=\Gamma\backslash\br{I}_{p,n-p+1}$. Notice that \eqref{eq 21.1} has level one, with $\underline{h}\big(\V^{\w_{1}}\big\{ \tfrac{a}{2}\big\} \big)=(p,n-p+1)$ and $\underline{h}\big(\V^{\w_{n}}\big\{ \tfrac{-a}{2}\big\} \big)=(n-p+1,p)$ for its two complex summands. (These are both the Hodge numbers and the signatures of the relevant Hermitian forms.) Our focus on these cases is motivated by the
\begin{Definition}[cf.~\cite{Iz}] \label{DefGW} A \emph{universal $k$-Weil abelian $(n+1)$-fold} is an abelian variety $\CA\to\mathrm{X}$ whose $H^{1}$ recovers $\H_{p,n-p+1}$ with $n-p+1=p+k$.
\end{Definition}

As in the ($0$-)Weil case, $\operatorname{End}(\CA)_{\QQ}$ contains an imaginary quadratic field; but the generic Mumford--Tate group is $U(p,p+k)$ instead of ${\rm SU}(p,p)$. For the $H_{{\rm prim}}^{2d+1}$ of~$\CA$, an easy argument with Young diagrams shows that\footnote{Where $\w_{0}+\w_{n-2d}$ [resp.~$\w_{2d+1}+\w_{0}$] means $\w_{n-2d}$ [resp.~$\w_{2d+1}$].}
\begin{gather} \label{eq 21.2}
\frac{\bigwedge^{2d+1}\H_{p,n-p+1}}{\bigwedge^{2d-1}\H_{p,d-p+1}} \cong \bigoplus_{j=0}^{2d+1} \V^{\w_j+\w_{n-2d+j}}\left\{\frac{(n-2p+1)(2d-2j+1)}{2(n+1)} \right\}
\end{gather}for $1\leq d\leq\tfrac{n}{2}$. There are no terms of level one in \eqref{eq 21.2}, except at $j=0,2d+1$ when $2d+1=n+1$ in case $\br{I}_{p,p+1}$ and when $2d+1=n$ in case $\br{I}_{p,p}$. For all the other terms on the right-hand side of \eqref{eq 21.2}, we compute using \eqref{eq 15.1}
\begin{gather*}
\mu \left( \w_j + \w_{n-2d+j} , \frac{(n-2p+1)(2n-2j+1)}{n+1} \right) = \begin{cases} j-d-p, & p<j \\ 1-d, & p=j\text{ or }n-2d+j, \\ -d, & j<p<n-2d+j, \\ p+d-j-n, & n-2d+j<p. \end{cases} 
\end{gather*}Since $p+d-j-n\geq0$ (with $d\leq\tfrac{n}{2}$, $p\leq\tfrac{n+1}{2}$) boils down to the level one cases we identified, the last two entries can be ignored. We have $1-d\geq0$ and $p=j$ or $n-2d+j$ in the cases $(d,p)=(1,1)$, $(1,2)$, $(1,3)$; and $j-d-p\geq0$ for some $j$ $\implies$ $2d+1\geq2p-1$. Hence we conclude that the variations of level~$>1$ in $H^{{\rm odd}}$ of $\CA$ admitting infinitesimal normal functions therefore lie in
\begin{gather*}
 \begin{cases}
\text{all possible degrees }(3\text{ through }2n-1), & \text{if }p\leq3,\\
\text{degrees }2p-1\text{ through }2n-2p+3, & \text{if }3<p<\tfrac{n+1}{2},\\
\text{no degrees}, & \text{if }3<p=\tfrac{n+1}{2}.
\end{cases}
\end{gather*}
Note that $p=\tfrac{n+1}{2}$ ($k=0$) is the Weil abelian setting dispensed with in Section~\ref{S3.5}. Omitting this case, we can restate our
conclusion as follows:
\begin{Theorem}\label{thm3}For $k\leq n-7$, no \'etale pullback of a universal $k$-Weil abelian $(n+1)$-fold admits reduced normal functions arising from cycles of dimension less than $\tfrac{n-k-1}{2}$ or codimension less than $\tfrac{n-k+1}{2}$.
\end{Theorem}
As to why (or whether) greater disparity in the signature of the Hermitian form should lead to more $\overline{\operatorname{AJ}}$-nontrivial cycles in the Griffiths group, we can say nothing yet.

\subsection{Comments on the Hermitian Calabi--Yau VHS} \label{S4.2}

For the tube domain cases, the CY variations are nothing but the $\V_{\RR}^{k\w_{{\tt I}}}$ (described for $k=1$ by Gross~\cite{Gro}). In the remaining cases (after \cite{FL, SZ}), where $V^{k\w_{{\tt I}}}$ is complex, we choose the shift in~$\tilde{\V}^{k\w_{{\tt I}}}\big\{ \tfrac{a}{2}\big\} $ to minimize the (odd) level while having the CY property.

Writing $\epsilon(m)=1$ ($m$ even) resp.~$2$ ($m$ odd), to accomplish this we need (cf.~\eqref{eq 6.3} and \eqref{eq 20.2})
\begin{gather*}
\ell\big(\tilde{V}^{k\w_{{\tt I}}}\big\{ \tfrac{a}{2}\big\} \big)=\ell\big(V^{k\w_{{\tt I}}}\big)+\epsilon\big(\ell\big(V^{k\w_{{\tt I}}}\big)\big)\ \implies\\
a=\tfrac{1}{2}\big\{ {\tt E}(k\w_{{\tt I}})-{\tt E}(k\tau(\w_{{\tt I}}))\big\} +\epsilon\big(\tfrac{1}{2}\big\{ {\tt E}(k\w_{{\tt I}})+{\tt E}(k\tau(\w_{{\tt I}}))\big\} \big).
\end{gather*}
This yields $a=\tfrac{2kp^{2}}{n+1}-kp+\epsilon(kp)$ (for $\br{I}_{p,n-p+1(\neq p)}$), $a=-\tfrac{k}{2}+\epsilon(mk)$ (for $\br{II}_{2m+1}$), and $a=-\tfrac{2k}{3}+1$ (for $\br{EIII}$).\footnote{For $k=1$, this gives the complex examples in~\cite{FL}, where we note that the sign on our half-twist is the opposite of theirs. Note that [op.\ cit.] omits the $(p,k)=(1,2)$ cases $\br{I}_{1,n} (A_{n},\sigma_{1};2\w_{1} )\big\{ \tfrac{3-n}{2(n+1)}\big\}$, which also have level~3.} Plugging into \eqref{eq 20.3} gives
\begin{gather*}
\mu(k\w_{{\tt I}},a)=k-\big\lfloor \tfrac{kp+1}{2}\big\rfloor ,\qquad 1-\big\lfloor \tfrac{mk+1}{2}\big\rfloor ,\qquad \text{resp.}\ 1-k,
\end{gather*}
which is non-negative when $p=1$ or $2$ ($\br{I}_{p,n-p+1}$), $(m,k)=(2,1)$ ($\br{II}_{2m+1}$), resp.~$k=1$ ($\br{EIII}$). This includes all the level 3 cases, but also two infinite series of examples (since $k$ is arbitrary) for $\br{I}_{1,n}$ and $\br{I}_{2,n-1}$. Putting this together with the results in Section~\ref{S3}, we have the
\begin{Proposition} Amongst the minimal-level CY variations over irreducible Hermitian symmetric domains other than $\br{I}_{1,n}$ or $\br{I}_{2,n-1}$, only those of level three admit infinitesimal normal functions.
\end{Proposition}

\subsection*{Acknowledgments}

The authors thank P.~Brosnan and G.~Pearlstein for helpful discussions, the referees for their careful reading, and gratefully acknowledge support from NSF Grant DMS-1361147. This paper was written while MK was a member at the Institute for Advanced Study, and he thanks the IAS for excellent working conditions and the Fund for Mathematics for financial support.

\pdfbookmark[1]{References}{ref}
\LastPageEnding

\end{document}